\documentclass[a4paper]{jpconf}
\usepackage{graphicx}
\begin{document}
\title{ Hypersurfaces with many $A_{j}-$singularities: explicit constructions.}
 \author{Juan Garc\'{\i}a Escudero}

\address{Universidad de Oviedo. Facultad de Ciencias,
 33007 Oviedo, Spain}

\begin{abstract}
A construction of algebraic surfaces based on two types of simple arrangements of lines, containing the prototiles of substitution tilings, has been proposed recently. The surfaces are derived with the help of polynomials obtained from products of the lines generating the simple arrangements. One of the arrangements gives the generalizations of the Chebyshev polynomials known as folding polynomials. The other produces a family of polynomials having more critical points with the same critical values, which can also be used to derive hypersurfaces with many $A_{j}-$singularities.  
\end{abstract}

%

 \section {Introduction.}  
  \bigskip\par  
  Algebraic hypersurfaces with many $A_{j}$-singularities have been studied in \cite{lab06} by extending a construction of S.V. Chmutov \cite{chm92}, who used generalizations of the Chebyshev polynomials known as folding polynomials \cite{hof88,wit88}.  Chmutov constructions produce surfaces having a high number of nodes or $A_{1}$-singularities, although for degrees $d=6,7,8,10,12$ there are examples with a higher number \cite{bar96,end97,lab06b,sar01}. Nodal surfaces of several degrees have also been used  to construct  two-dimensional Calabi-Yau manifolds, also known as K3-surfaces. This is the case for the 600 nodal Sarti dodecic, a surface having a quotient which is a K3-surface \cite{bar03}. New lower bounds for the number of singularities $A_{j}$ on a hypersurface of degree $d$ in the complex projective space ${\bf{P}}^{n}({\bf{C}})$ are obtained for many cases in  \cite{lab06}.
  \par
Two types of simple arrangements  \cite{gru72} of lines $\Sigma_{1}, \Sigma_{2}$, containing the prototiles of substitution tilings, have been introduced in \cite{esc11c}. Real variants of Chmutov surfaces \cite{bre08} can be derived with the folding polynomials defined as products of lines in $\Sigma_{1}$. The mirror symmetric simple arrangements $\Sigma_{2}$,  $\bar{\Sigma}_{2}$, with cyclic symmetries $C_{3}$, have one more triangle than $\Sigma_{1}$, and this property can be used to construct surfaces with more singularities. This is due to the fact that the corresponding polynomials have the same extreme values (-1 for all the minima and 8 for the maxima) in all the critical points with critical value of the same sign, which on the other hand coincide with those of the real folding polynomials constructed with $\Sigma_{1}$. The triangular shapes appearing in the simple arrangements are the prototiles of substitution tilings, derived with a general construction based on simplicial arrangements of lines \cite{esc08}. A wide variety of tiling spaces can be defined having different types of topological invariants \cite{esc11a,esc11b}, which are significant in several contexts, like the study of the energy gaps in the spectrum of the Schroedinger equation with tiling equivariant potentials. In this work, the construction of new types of hypersurfaces with many $A_{j}$-singularities is done with the polynomials associated with  $\Sigma_{2}$ and $\bar{\Sigma}_{2}$. Mathematica \cite{wol91}, Singular \cite{gre08}  and Surfer \cite{end01,hol05} computing and geometric visualization tools are used.
  
 \section {Surfaces with many $A_{1}$-singularities}  
  \bigskip\par  
 Let $T_{d}^{C}(w)\in {\bf{R}}[w]$ be the Chebyshev polynomial of degree $d$ with critical values -1,+1 and $F_{d}(u,v)$ the folding polynomials associated to the system of roots ${\bf{A}}_{2}$. An explicit formula for the folding polynomials of degree $d$, which are symmetric ($F_{d}(u,v)=F_{d}(v,u)$) with real coefficients, is $F_{d}(u,v)=2+ j_{d}(u,v)+ j_{d}(v,u)$, where
    \par
\[ j_{d}(u,v)=\det \left(
    \begin{array}{llllllll}
 u &1 &0 &0&0 &0 &... &0   \\
    2v &u &1 &0&0 &0 &... &0   \\
     3 &v &u &1 &0 &0 &... &0  \\
        0 &1 &v &u&1 &0 &... &0   \\
               0 &0 &1 &v&u &1 &0.. &0   \\
                      ... &... &... &...&... &... &... &0   \\
                             ... &... &... &0&1 &v &u &1   \\
                                    0 &... &... &0&0 &1 &v &u   \\

    \end{array}
    \right) \]
  \bigskip\par
  Chmutov constructs surfaces $Ch_{d}(u,v,w):=F_{d}(u,v)+(T_{d}^{C}(w)+1)/2$ having many nodes. In \cite{bre08} it is shown that the Chmutov construction can be adapted to give only real singularities. The real folding polynomials $F_{d}(x,y)$ are obtained when $u=x+iy,v=u^{\ast}=x-iy$. The plane curve $F_{d}(x,y)$ is a product of $d$ lines in simple arrangements $\Sigma_{1}$ having critical points with only three different critical values: 0,-1 and 8.  The surface $Ch_{d}(x,y,z)$ is singular exactly at the points where the critical values  $\zeta$ of $F_{d}(x,y)$ and $(T_{d}^{C}(z)+1)/2$ are either both zero, or one is -1 and the other +1.
\par
The polynomials $F_{d}(x,y)$ have ${d \choose 2}$ real critical points with value $\zeta=0$ and, when  $d \in {\bf{Z}}_{3}$, $\frac{1}{3}d^{2}-d$ real critical points with  $\zeta=-1$. The other critical points also have real coordinates and critical value $\zeta=8$. The number of maxima can then be obtained by applying the following
  \bigskip\par
  {\it{  Lemma 1.}}  \cite{ort03,bre08} Let $f$ be a real simple line arrangement consisting in $d \ge 3$ lines.  $f$  has exactly ${d-1 \choose 2}$ bounded open cells each of which contains exactly one critical point. All the critical points of $f$ are non-degenerate.
 \bigskip\par
It is possible to construct polynomials for some degrees with one more critical point with  $\zeta=-1$. 
For $m=3q, q=1,2,3,...$ the simple arrangement $\bar{\Sigma}_{2}$ consists in the lines $\bar{L}_{k,m}(x,y)=0$, with
     \begin{equation}
\bar{L}_{k,m}(x,y):=y-({\rm cos} \frac{k \pi}{3m}-x){\rm tan}\frac{k \pi}{6m}-{\rm sin} \frac{k  \pi}{3m}
 \end{equation}   
 where $k\in \bf{Z}$ and $\bar{L}_{3m,m}(x,y)=0$ is interpreted as the line $x=-1$. The polynomials obtained with $\bar{L}_{k,m}(x,y)$ are defined as 
      \begin{equation}
\bar{J}_{m}^{C}(x,y):= 3^{\frac{1-(-1)^{m}}{4}} (-1)^{\lfloor \frac{q+1}{2}\rfloor+1} \prod_{\nu=0}^{m-1}\bar{L}_{6\nu+1,m}(x,y)\in {\bf{R}}[x,y]
 \end{equation} 
and $J_{m}^{C}(x,y)=\bar{J}_{m}^{C}(x,-y)$. They have only three different critical values: 0,-1, 8, and by changing the variables $x=(u+v)/2, y=i(v-u)/2$ in $J_{m}^C(x,y)$,  we get
$$J_{m}^C(u,v)=-1+ j_{m}^{C}(u,v)+  j_{m}^{C \ast}(u,v)$$
\par\noindent
where $ j_{m}^{C}(u,v)$ has now complex coefficients. The polynomials corresponding to $d=3,6,9$, are  
$$
j_{3}^{C}(u,v)=-bu^3+\frac{3}{2}uv,$$

$$
j_{6}^{C}(u,v)=(2-8b)u^3-bu^6+6uv+6bu^{4}v-\frac{9}{2}u^{2}v^{2},$$

$$
j_{9}^{C}(u,v)=(9-27b)u^3-9bu^6-bu^9+\frac{27}{2}uv-(9-54b)u^{4}v+9bu^{7}v-27u^{2}v^{2}-27bu^{5}v^{2}+15u^{3}v^{3},$$

  \bigskip\par
\par\noindent 
with $b=e^{-i\frac{\pi}{3}}$.
In contrast with $F_{d}(u,v)$, the polynomials $J_{m}^C(v,u)$  are not symmetric and $\bar{J}_{m}^C(u,v):={J}_{m}^C(u^{\ast},v^{\ast})=J_{m}^C(v,u)$. The surfaces with affine equations \cite{esc11c}
  \begin{equation}
Q_{m}^{C}(x,y,z):=J_{m}^{C}(x,y)+\frac{(-1)^{m+1}}{4}(J_{m}^{C}(z,0)-1+(-1)^{m+1}2)=0
\end{equation}
 \par\noindent
have a number of real nodes, or $A_{1}$-singularities, higher than the real variants of the Chmutov surfaces with the same degree. The first surface in the series defined by Eq.(3) is equivalent to the Cayley cubic, a well known example having four ordinary double points,  the maximum possible number in degree 3. The Chmutov surface of the same degree has three nodes.

 \section {Hypersurfaces with many $A_{j}$-singularities.}  
  \bigskip\par  
An $A_{j}$-singularity on a hypersurface in ${\bf{P}}^{n}({\bf{C}})$ has the local equation $x_{0}^{j+1}+x_{1}^2+x_{2}^2+...+x_{n-1}^2=0$. The polynomials $J_{m}^{C}(x,y)$ and $\bar{J}_{m}^{C}(x,y)$ can be used for the construction of hypersurfaces with many $A_{j}$-singularities.

 \begin{figure}[h]
 \includegraphics[width=20pc]{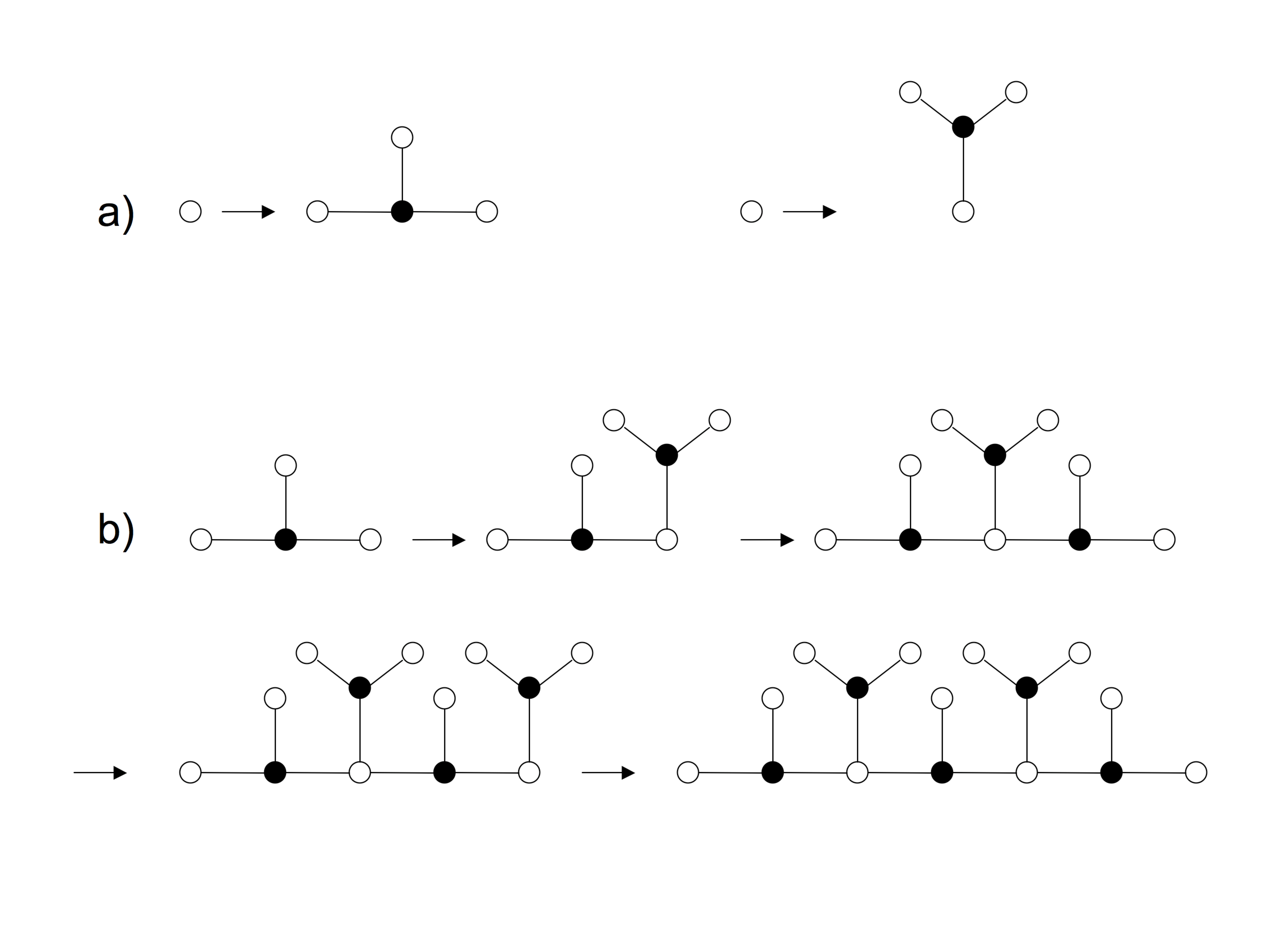}
\caption{\label{label} Plane trees for Belyi polynomials.(a)The two types of  substitution rules. (b) Plane trees corresponding to $\bar{B}_{3q}^{2}(z)$, $q=1,2,3,4,5$, obtained by applying alternatively  the  substitutions to the right most vertex.}
\end{figure}

  \par
  A polynomial in one variable with no more than two different critical values is called a Belyi polynomial. The proof of its existence is based on the theory of Dessins d\'{}Enfants. A graph without cycles with a prescribed cyclic order of the edges adjacent to each vertex is called a plane tree \cite{adr98}. A plane tree has a natural bicoloring of the vertices. For any plane tree, there exists a Belyi polynomial whose critical points have the multiplicities given by the number of edges adjacent to the vertices minus one and viceversa. Degree-$d$ Chebyshev polynomials are the Belyi polynomials whose plane tree has the form of the Dynkin diagram of the Lie algebra ${\bf{A}}_{d+1}$. The plane trees of a family of Belyi polynomials that we denote by $\bar{B}_{3q}^{2}(z)$ (see Appendix B), can be obtained by a substitution process as indicated in Fig.1. After the first step, which is given in the left part of Fig.1a, we apply alternatively to the right most white vertex,  the two substitution rules given in Fig.1a. The next four steps, corresponding to $q=2,3,4,5$ are shown in Fig.1b. The generalization of the Chmutov\'{}s construction given in \cite{lab06} for the study of singularities of type $A_{j}$ in ${\bf{P}}^{3}({\bf{C}})$ consists in considering surfaces with affine equations 
  \begin{equation}   
    F_{d}(x,y)+G_{d}^{j}(z)=0
    \end{equation}
     where  $G_{d}^{j}(z)$ are Belyi polynomials.  In general there is no explicit formula for them, and they can be computed only for low degree $d$ by using Groebner basis. However, in some cases, there are explicit expressions that can be obtained from classical Jacobi polynomials $P_{k}^{l,m}(z)$ \cite{rob07} by studying the three-parameter family of polynomials 
   \begin{equation}    
     G_{a,b,c}(z)=z^{a}P_{b-1}^{a/c,-b}(1-2z)^{c}
   \end{equation}  
     
   The order of a zero $z_{0}$ of $\frac{dP(z)}{dz}$, which is a critical point of $P(z)$ with critical value $\zeta=P(z_{0})$, is called its multipicity $\nu$ (all the derivatives of $P(z)$ up to order $\nu$ vanish at $z_{0}$). The zeros of $\frac{dG_{a,b,c}}{dz}$ with critical value $\zeta=0$ have orders $\nu=a-1,c-1,...,c-1$, while the unique remaining root $r=1$ has $\nu=b-1$ with critical value $\zeta=1$. 
    \par
   We denote by $\mu_{A_{j}}$ the number of $A_{j}-$singularities ($\mu_{A_{j}}(d)$ if we want to specify the degree $d$).  If we replace $F_{3q}(x,y)$ by $J_{3q}^{C}(x,y)$ or  $\bar{J}_{3q}^{C}(x,y)$  in Eq.(4) we can get in some cases surfaces with a higher $\mu_{A_{j}}$. In what follows, we use $G_{a,b,c}(z)$ in order to give explicit expressions for some singular surfaces. Even when the number of singularities of a given type is equal to the one obtained in \cite{lab06}, the surfaces show a higher number of other types. The lines $\bar{L}_{k,m}(x,y)=0$ and $\bar{L}_{l,m}(x,y)=0$, where
 
  \begin{equation}
 \bar{L}_{k,m}(x,y)=y+t_{k,m}x-B_{k,m}
 \end{equation}
   \par\noindent
 with $t_{k,m}={\rm tan}\frac{k\pi}{6m}, B_{k,m}=(3t_{k,m}-t_{k,m}^3)/(1+t_{k,m}^2)$, intersect in the point 
 \begin{equation}
(\frac{ B_{k,m}-B_{l,m}}{t_{k,m}-t_{l,m}},\frac{ B_{l,m}t_{k,m}-B_{k,m}t_{l,m}}{t_{k,m}-t_{l,m}})
\end{equation}
which is denoted, for a given $m$, by $p_{k\cap l}$. 
  \bigskip\par  
  {\it{  Proposition 1.}}  The surfaces $J_{6}^{C}(x,y)+G_{2,2,4}(z)=0$ in ${\bf{P}}^{3}({\bf{R}})$ have $\mu_{A_{1}}=22$ and $\mu_{A_{3}}=15$.  In ${\bf{P}}^{4}({\bf{R}})$ the three-fold $J_{6}^{C}(x_{0},x_{1})-J_{6}^{C}(x_{2},x_{3})=0$ has $\mu_{A_{1}}=283$.
  \bigskip\par  
  {\it{  Proof.}}     
  The points corresponding to the minima of $J_{6}^{C}(x,y)$ are, up to three-fold rotation, $p_{0\cap 6},p_{0\cap 12},p_{2\cap 8}$ with critical value $\zeta=-1$ and the maximum with critical value $\zeta=8$ is in $p_{4\cap 10}$.  Except for the minima  situated at the origin ($p_{0\cap 12}$ in this case), the action of the cyclic group $C_{3}$ gives the remaining distinct critical points for this and the examples that follows. By representing the points in the plane by complex numbers, we find that the minima are in
  $$0,\alpha, r \alpha , r^{2} \alpha, 2, 2r, 2r^{2} $$
  where $\alpha=e^{i\frac{\pi}{9}},r=e^{i\frac{2 \pi}{3}}$. The maxima correspond to the conjugates of the minima which are not themselves minima. In this case the maxima are in
  $$\alpha^{\ast}, r^{\ast} \alpha^{\ast} , r^{\ast 2} \alpha^{\ast}$$
\par\noindent
and $Q_{6}^{C}(x,y,z)$ has Tjurina number $\tau=59$ (see Appendix A). The polynomial $G_{2,2,4}(z)=\frac{{\left( -3 + z \right)
        }^4\,z^2}{16}$  has critical points $z=0,3$ (multiplicity $\nu=3$) with $\zeta=0$ and $z=1$ with $\zeta=1$. Although 15 is the best known lower bound of $\mu_{A_{3}}$ for a sextic in ${\bf{P}}^{3}({\bf{C}})$,  the surface presented here has higher  $\mu_{A_{1}}$ and the singularities are real (Fig.2a). 
    \bigskip\par  
  {\it{  Proposition 2.}} In  ${\bf {P}}^{4}({\bf{R}})$ the hypersurface $J_{9}^{C}(x_{0},x_{1})-J_{9}^{C}(x_{2},x_{3})=0$ has $\mu_{A_{1}}=1738$. The following polynomials define surfaces of type $J_{9}^{C}(x,y)+G_{a,b,c}(z)=0$ in ${\bf{P}}^{3}({\bf{K}})$ with the indicated number of singularities: 
    \par 
2.1.- $G_{3,3,3}(z)$; ${\bf{K}}={\bf{C}}$ ;  $\mu_{A_{2}}=127$. 
     \par 
2.2.- $G_{1,3,4}(z)$; ${\bf{K}}={\bf{C}}$;  $\mu_{A_{3}}=72, \mu_{A_{2}}=19$ .    
\par 
2.3.- $G_{3,2,6}(z)$; ${\bf{K}}={\bf{R}}$;  $\mu_{A_{5}}=\mu_{A_{2}}=36,\mu_{A_{1}}=19$.
    \par 
The maximum number of cusps and real $A_{4}$-singularities for a nonic surface satisfy $\mu_{A_{2}}(9)\ge 127, \mu_{A_{4}}(9)\ge 55$.
  \bigskip\par  
  {\it{  Proof.}}   
  The result for the hypersurface in ${\bf {P}}^{4}({\bf{R}})$ is obtained having in mind that the minima of  $J_{9}^{C}(x,y)$ are in $p_{0\cap 6},p_{0\cap 12},p_{0\cap 18},p_{0\cap 24},p_{2\cap 8},p_{2\cap 14},p_{8\cap 14}$ with $\zeta=-1$ and the maxima with $\zeta=8$ are situated in $p_{4\cap 10},p_{4\cap 16},p_{10\cap 16}$. The number of singularities of the surfaces in ${\bf{P}}^{3}({\bf{K}})$ can be computed by using the following properties:
   \par
1.- $G_{3,3,3}(z)=      z^3\,{\left( 3 - 3\,z + 
      z^2 \right) }^3$  has critical points $z=0,\frac{3 \pm i \,{\sqrt{3}}}
  {2}$ with $\zeta=0$ and $z=1$ with $\zeta=1$, all of them with $\nu=2$.   
   \par
2.- $G_{1,3,4}(z)=  z\,{\left( \frac{45}{32} - 
      \frac{9\,z}{16} + 
      \frac{5\,z^2}{32}
      \right) }^4$  has critical points $z=\frac{9 \pm 12 i}
  {5}$  ($\nu=3$) with $\zeta= 0$ and $z=1$ ($\nu=2$) with $\zeta=1$.
   \par  
3.- $G_{3,2,6}(z)= \frac{{\left( -3 + z \right)}^6\,z^3}{64}$  has critical points $z=0$ ($\nu=2$), $3$ ($\nu=5$) with $\zeta=0$ and $z=1$ with $\zeta=1$.  The surface $J_{9}^{C}(x,y)+G_{3,2,6}(z)=0$ is represented in Fig.2b.
  \par
By using other Belyi polynomials instead of $G_{3,3,3}(z)$, one gets surfaces with $\mu_{A_{2}}=110$ (see Appendix B). The existence of the surface $J_{9}^{C}(x,y)+G_{3,3,3}(z)=0$ increases in one cusp the known lower bound of $\mu_{A_{2}}$ for a nonic \cite{lab06}. The Belyi polynomial $\bar{B}_{9}^{4}(z)$ given in Appendix B produces a nonic with $\mu_{A_{4}}=55$.

    \bigskip\par  
  {\it{Proposition 3.}} The three-fold $J_{12}^{C}(x_{0},x_{1})-J_{12}^{C}(x_{2},x_{3})=0$ in ${\bf{P}}^{4}({\bf{R}})$ has $\mu_{A_{1}}=6049$. The number of singularities of $J_{12}^{C}(x,y)+G_{a,b,c}(z)=0$ in ${\bf{P}}^{3}({\bf{K}})$ is  
    \par 
3.1.-   $\mu_{A_{4}}=132,\mu_{A_{2}}=37,\mu_{A_{1}}=66$ for $a=2,b=3,c=5$; ${\bf{K}}={\bf{C}}$.
     \par 
3.2.-  $\mu_{A_{7}}=\mu_{A_{3}}=66, \mu_{A_{1}}=37$ for $a=4,b=2,c=8$; ${\bf{K}}={\bf{R}}$.  
\par
The following lower bound is valid for dodecic surfaces:  $\mu_{A_{2}}(12)\ge 301$. 
  \bigskip\par  
  {\it{  Proof.}}  For degree-12 we have  
  
$$
j_{12}^{C}(u,v)=(24-64b)u^3-(2+40b)u^6-12bu^9-bu^{12}+24uv-(60-240b)u^{4}v+96bu^{7}v+12bu^{10}v$$
$$-90u^{2}v^{2}-(36-288b)u^{5}v^{2}-54bu^{8}v^{2}+120u^{3}v^{3}-28bu^{6}v^{3}-\frac{105}{2}u^{4}v^{4},$$ 
with $b=e^{-i\frac{\pi}{3}}$. We denote by $\it{S}[A]$ the 2-combinations or subsets of two elements $\{k,l\}$ taken from the set $\it{A}$. The values $\{k,l\}$ in $p_{k\cap l}$ for the minima of  $J_{12}^{C}(x,y)$ with $\zeta=-1$ are $\{6,12\},\{66,60\}$,$\{0,6n\},n=1,2,3,4,5$ together with the ones specified by the elements of ${\it{S}}[\{2,8,14,20\}]$. The maxima with $\zeta=8$ are obtained from ${\it{S}}[\{4,10,16,22\}]$. Now we have:
 \par
1.-$G_{2,3,5}(z)= z^2\,{\left( \frac{42}{25} - 
      \frac{24\,z}{25} + 
      \frac{7\,z^2}{25}
      \right) }^5$  has critical points $z=0,\frac{12 \pm i 5\,{\sqrt{6}}}
  {7}$ ($\nu=4$) with $\zeta=0$ and $z=1$ ($\nu=2$) with $\zeta=1$. 
   \par
2.-$G_{4,2,8}(z)= \frac{{\left( -3 + z \right)
        }^8\,z^4}{256}$  has critical points $z=0$ ($\nu=3$),$3$ ($\nu=7$) with $\zeta=0$ and $z=1$ with $\zeta=1$. 
         \par
  The existence of $J_{12}^{C}(x,y)$ and the Belyi polynomial with the plane tree as in Fig.1b, gives a surface with $\mu_{A_{2}}(12)=301$.
    \bigskip\par  
  {\it{Proposition 4.}} In ${\bf{P}}^{4}({\bf{R}})$  the three-fold $J_{15}^{C}(x_{0},x_{1})-J_{15}^{C}(x_{2},x_{3})=0$ has $\mu_{A_{1}}=15646$. The number of singularities of $J_{15}^{C}(x,y)+G_{a,b,c}(z)=0$ in ${\bf{P}}^{3}({\bf{K}})$ is    
  \par 
4.1.- $\mu_{A_{3}}=376,\mu_{A_{2}}=61$ for $a=3,b=4,c=4$; ${\bf{K}}={\bf{C}}$.
     \par 
4.2.-  $\mu_{A_{5}}=210, \mu_{A_{2}}=166$ for $a=3,b=3,c=6$; ${\bf{K}}={\bf{C}}$.    
    \par 
4.3.-  $\mu_{A_{9}}=\mu_{A_{4}}=315, \mu_{A_{1}}=61$ for $a=5,b=2,c=10$; ${\bf{K}}={\bf{R}}$.
\par
We have the following lower bounds for 15-degree surfaces:  $\mu_{A_{2}}(15)\ge 647, \mu_{A_{3}}(15)\ge 376$. 
  \bigskip\par  
  {\it{  Proof.}}   In this case  
$$j_{15}^{C}(u,v)=(50-125b)u^3-(15+125b)u^6-75bu^9-15bu^{12}-bu^{15}+\frac{75}{2}uv-(225-750b)u^{4}v$$ 
$$+(15+525b)u^{7}v+165bu^{10}v+15bu^{13}v-225u^{2}v^{2}+(315-1575b)u^{5}v^{2}-675bu^{8}v^{2}-90bu^{11}v^{2}$$
$$+525u^{3}v^{3}-(140-1400b)u^{6}v^{3}+275bu^{9}v^{3}-525u^{4}v^{4}-450bu^{7}v^{4}+189u^{5}v^{5},$$
 with $b=e^{-i\frac{\pi}{3}}$.  The subsets $\{k,l\}$  associated to the minima of  $J_{15}^{C}(x,y)$ with $\zeta=-1$ are $\{6,12\}, \{6,18\},\{84,78\},\{84,72\}$,$\{0,6n\},n=1,2,...,7$,  as well as those given by the elements of ${\it{S}}[\{2,8,14,20,26\}]$. The maxima with $\zeta=8$ are obtained from ${\it{S}}[\{4,10,16,22,28\}]$. The polynomials  $G_{a,b,c}(z)$ considered here have the following properties:
  \par    
1.-$G_{3,4,4}(z)= z^3\,{\left( \frac{385}
       {128} - 
      \frac{495\,z}{128} + 
      \frac{315\,z^2}{128} - 
      \frac{77\,z^3}{128}
      \right) }^4$  has critical points $z=0$ ($\nu=2$), $z\approx 2.20309$ ($\nu=3$), $z\approx 0.94391\pm i 1.17413$ ($\nu=3$) with $\zeta=0$ and $z=1$ ($\nu=3$) with $\zeta=1$. 
 \par       
          
2.-$G_{3,3,6}(z)= z^3\,{\left( \frac{15}{8} - 
      \frac{5\,z}{4} + 
      \frac{3\,z^2}{8} \right)
      }^6$  has critical points $z=0$ ($\nu=2$), $z=\frac{3 \pm i 2\,{\sqrt{5}}}
  {7}$ ($\nu=5$) with $\zeta=0$ and $z=1$ ($\nu=2$) with $\zeta=1$. 
\par    

3.-$G_{5,2,10}(z)=\frac{{\left( -3 + z \right)
        }^{10}\,z^5}{1024}$  has critical points $z=0$ ($\nu=4$), $3$ ($\nu=9$) with $\zeta=0$ and $z=1$ with $\zeta=1$. 
        \par
  The Belyi polynomial $\bar{B}_{15}^{2}(z)$ associated with the plane tree in Fig.1b, together with $J_{15}^{C}(x,y)$ gives a surface with $\mu_{A_{2}}(15)=647$.  We notice that $J_{15}^{C}(x,y)+G_{3,4,4}(z)=0$ has one more $A_{3}$-singularity than the surfaces with the same degree having the highest known $\mu_{A_{3}}$ \cite{lab06}.     
Continuing along this lines, degree-$3q$ surfaces in ${\bf{P}}^{3}({\bf{R}})$ with singularities of type $A_{2q-1}, q=1,2,...$ can be obtained by adding $G_{q,2,2q}(z)$ to  $J_{m}^{C}(x,y)$ or $\bar{J}_{m}^{C}(x,y)$.
   \bigskip\par 
 Nodal hypersurfaces in ${\bf{P}}^{n}({\bf{R}})$ can be constructed also with the help of $J_{m}^{C}(x,y)$:        
 
  \begin{equation}
 \sum_{k=0}^{\lfloor\frac{n-2}{2}\rfloor}(-1)^{k}J_{m}^{C}(x_{2k},x_{2k+1})=\frac{(-1)^{m+n}+(-1)^{m+1}}{8}(J_{m}^{C}(x_{n-1},0)-1+(-1)^{m+1}2)
\end{equation}
  \par\noindent
For instance $J_{6}^{C}(x_{0},x_{1})-J_{6}^{C}(x_{2},x_{3})=\frac{1}{4}(3-J_{6}^{C}(x_{4},0))$ in ${\bf{P}}^{5}({\bf{R}})$ has 1059 nodes. These can be computed having in mind the results in Prop.1 and the fact that
$\frac{1}{4}(3-J_{6}^{C}(x_{4},0))$ has 3 critical points with $\zeta=0$ and 2 critical points with $\zeta=1$.  
\par
The methods given above can be extended in order to construct hypersurfaces in ${\bf{P}}^{n}({\bf{C}}), n>4$ with many $A_{j}$-singularities, by using $G_{a,b,c}(x_{n-1})$ instead of the right-hand side of Eq.(8). In this way we get, for example, the hypersurface $J_{6}^{C}(x_{0},x_{1})-J_{6}^{C}(x_{2},x_{3})=G_{2,2,4}(x_{4})$ in ${\bf{P}}^{5}({\bf{R}})$ with 388 nodes and 283 $A_{3}$-singularities.
\par
 There is a dynamical formulation in terms of a uniparametric family of line configurations \cite{esc11b} which can be used to get deformations of the surfaces, where some singularities disappear and others arise. It gives a way to transform the surfaces based on $\Sigma_{1}$ into the surfaces based on $\Sigma_{2}$ and $\bar{\Sigma}_{2}$, representing a topology change that would correspond to a kind of phase transition in a physical context.
 \par 
We have studied hypersurfaces with a high number of singularities, improving the known lower bounds in some cases. Open questions that should be addressed are the search of maximal simple arrangements of lines giving polynomials with the same critical values, and deformations of the surfaces to get others with new types of singularities.

\bigskip\par 

 \begin{figure}[h]
 \includegraphics[width=35pc]{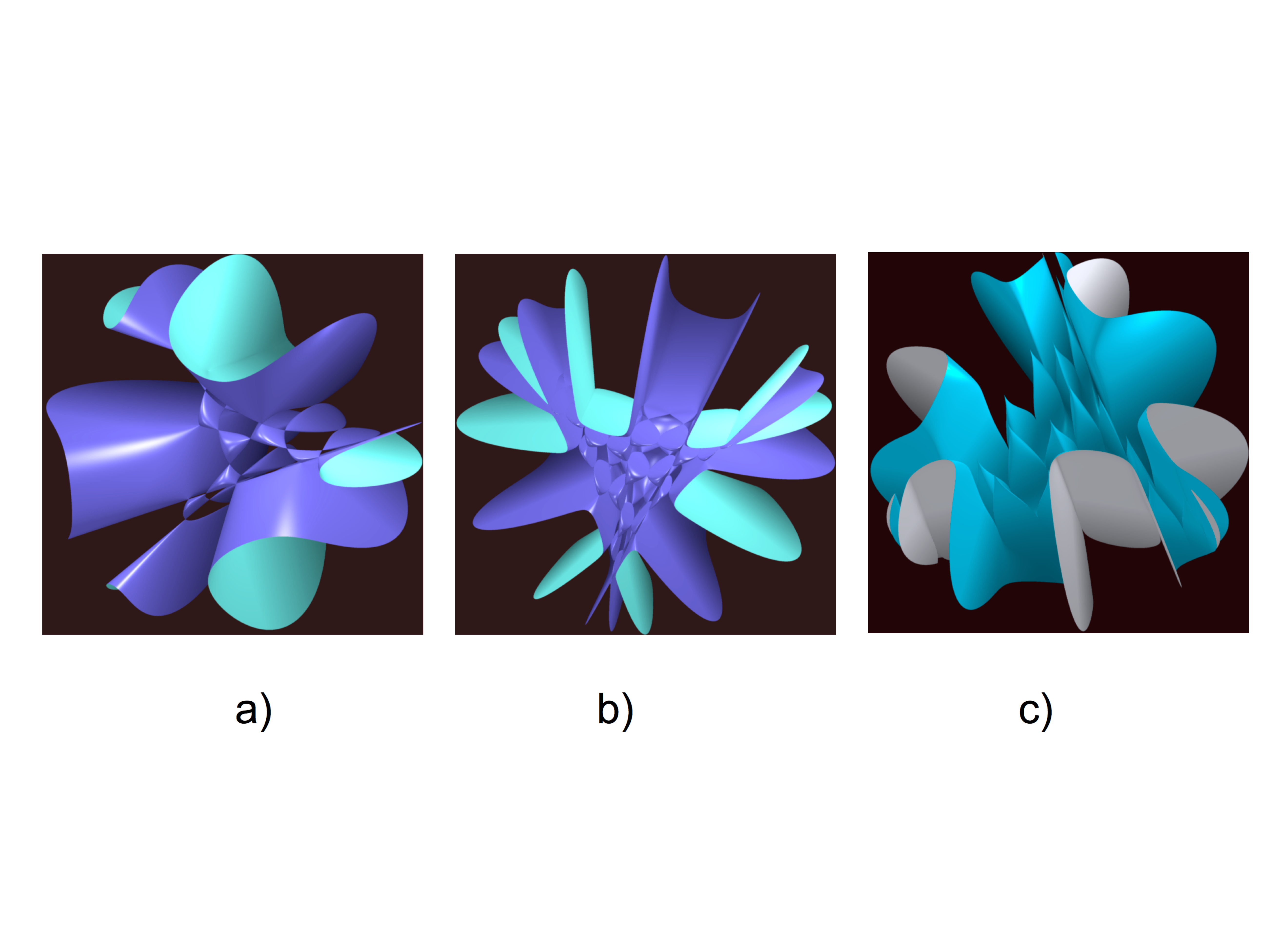}
\caption{\label{label} Examples with real singularities: a) A degree 6 surface with $\mu_{A_{1}}=22$ and $\mu_{A_{3}}=15$.b) A nonic surface with $\mu_{A_{1}}=19$ and $\mu_{A_{2}}=\mu_{A_{5}}=36$. c) A nonic with $\mu_{A_{4}}=55$.}
\end{figure}

 \section {APPENDIX A: Milnor and Tjurina numbers.}
 \bigskip\par 
Let $f$ be a holomorphic complex function germ at a given point. By $O$ we denote the ring of function germs and by 
$${\rm{J}}_{f}:=<\frac{ \partial f}{ \partial x_{1}},  \frac{ \partial f}{ \partial x_{2}},... \frac{ \partial f}{ \partial x_{n}}>$$
\noindent
 the jacobian ideal of $f$. The Milnor algebra $M$ of $f$ is given by the quotient algebra $O/\rm{J}_{f}$ and the Milnor number is given by the complex dimension of $M$. The Tjurina number $\tau$ is the dimension of the algebra obtained by replacing $\rm{J}_{f}$ by the ring  
 $$<f, \frac{ \partial f}{ \partial x_{1}},  \frac{ \partial f}{ \partial x_{2}},... \frac{ \partial f}{ \partial x_{n}}>$$.
\par   
 In what follows we use SINGULAR computer algebra system \cite{gre07,gre08} for checking the Tjurina and Milnor numbers.  The basering is a polynomial ring with three variables over the algebraic number field ${\bf{Q}}(a)$, $a=\sqrt{3}$. Lexicographical and degree reverse lexicographical monomial ordering are denoted by lp, dp, and
 vdim(std(I)) is the vector space dimension of the ring modulo de ideal I. Short input is used (e.g. $3x^2-x^3$ is denoted by $3x2-x3$). The result of the code
   \bigskip\par  
   LIB "sing.lib";
   \par  
ring R = (0,a),(x,y,z),dp;  //affine ring, , 0 is the characteristic of the ground field
\par  
minpoly=a2-3;//the minimal polynomial of a
\par  

poly f = -1+3x2-x3-3ax2y+3y2+3xy2+ay3+(3z2-z3)/4;

\par 

ideal sl = jacob(f),f; //the singular locus 
\par  
vdim(std(sl));    //Total Tjurina number 
   
  \bigskip\par\noindent  
is $\tau=4$, which is the number of singularities of the cubic $Q_{3}^{C}(x,y,z)=0$. For $Q_{d}^{C}(x,y,z), d=6,9,12,15$ we get $\tau=59,220,581,1162$.
\par 
  We can also use the built-in commands milnor and tjurina. For instance we can check that all the critical points of $J_{9}^{C}(x,y)$ are non degenerated (Lema 1): 
         \bigskip\par  
ring R = (0,alpha),(x,y),dp; 
  \par  
minpoly=a2-3;
  \par  
poly f = -1+27x2-9x3-54x4+36x5+21x6-27x7+9x8-x9-
81ax2y+162ax4y-54ax5y-81ax6y
\par 
+54ax7y-9ax8y
+27y2+27xy2-108x2y2-72x3y2+225x4y2+27x5y2-126x6y2+36x7y2
\par 
+27ay3+108ax2y3+180ax3y3-135ax4y3-126ax5y3+84ax6y3-
54y4-108xy4-45x2y4
\par 
+135x3y4-126x5y4
-54ay5-54axy5-27ax2y5-126ax3y5-126ax4y5
+39y6+81xy6+126x2y6
\par 
+84x3y6+
27ay7+54axy7+36ax2y7
-9y8-9xy8
-ay9;
  \par  
milnor(f);  
 \par  
ideal j = jacob(f);  // the critical locus of f
  \par  
poly h = det(jacob(j)); //determinant of the Hessian of f
  \par  
ideal nn = j,h; 
  \par  
vdim(std(nn)); 
 
  \bigskip\par\noindent
  gives a Milnor number of 64, and vdim(std(nn))=0.  
 For the sextic with non-nodes, the code
        \bigskip\par  

ring R = (0,a),(x,y,z),dp;  
\par  
minpoly=a2-3;
\par  
poly f =-1+12x2-4x3-9x4+6x5-x6-24ax2y+18ax4y-6ax5y+12y2+12xy2-18x2y2
\par  
-12x3y2+15x4y2+8ay3+12ax2y3+20ax3y3-9y4-18xy4-15x2y4-6ay5-6axy5+y6
\par 
+(z6-12z5+54z4-108z3+81z2)/16;
\par  
tjurina(f);  
  \bigskip\par\noindent
 gives $\tau=67$, as expected from Prop.1. In a similar way, if we analyze $f=J_{9}^{C}(x,y)+G_{3,3,3}(z)$, we get $\tau =254$, which corresponds to a surface with 127 cusps as in Prop.2.
 \section {APPENDIX B: Belyi polynomials.}
 \bigskip\par 
In this appendix we obtain several Belyi polynomials. Related with $A_{2}$-singularities for nonics are $B_{9}^{2}(z)$ and $\bar{B}_{9}^{2}(z)$.
       The roots of $\frac{dB_{9}^{2}(z)}{dz}$ with critical value  $\zeta=-1$ are denoted by $a,b$, while $0,u$ have critical value $\zeta=1$, all of them with multiplicity $\nu=2$. We have $\frac{dB_{9}^{2}(z)}{dz}=(z-a)^{2}(z-b)^{2}(z-u)^{2}z^{2}$, with $B_{9}^{2}(0)=1$. In order to get the values of the critical points we use the following code
      
  \bigskip\par 
ring R = 0,(a,b,u),lp;  
\par  
poly f1 = 2520+5a9-15a8b+12a7b2-15a8u+48a7bu-42a6b2u+12a7u2-42a6bu2+42a5b2u2;
\par  
poly f2 = 2520+12a2b7-15ab8+5b9-42a2b6u+48ab7u-15b8u+42a2b5u2-42ab6u2+12b7u2;
\par  
poly f3 = 42a2b2-42a2bu-42ab2u+12a2u2+48abu2+12b2u2-15au3-15bu3+5u4;
\par  
ideal I = f1, f2, f3;
\par  
ideal GI = groebner(I);
\par  
GI;
  \bigskip\par\noindent  

The Groebner basis has four elements: $GI_{1}=GI_{1}(u),GI_{2}=GI_{2}(u,b),GI_{3}=GI_{3}(u,b),GI_{4}=GI_{4}(u,b,a)$. 
 The polynomials $GI_{1},GI_{2}$, with degrees 180,172, have the common factor $u^{36}+135413275557888$. If we take one of its roots, for instance $u=2^{\frac{5}{9}} 3^{\frac{17}{36}} e^{\frac{i  \pi }{36}}$, then $a,b$ are the two complex roots of           
         $$ z^2 - u{2^{-\frac{3}{2}}}( 2^{\frac{3}{2}}  - 
       3^{ - \frac{1}{4}} - 
       3^{\frac{1}{4}} + 
  i (3^{- \frac{1}{4} } - 
          3^{\frac{1}{4}})) z
      + \frac{u^{2}}{24} ( 6 - 
        2^{\frac{1}{2}} 3^{\frac{5}{4}} - 
       2^{\frac{1}{2}} 3^{\frac{3}{4}} + 
      i (-2^{\frac{1}{2}} 3^{\frac{5}{4}} - 3^{\frac{1}{2}}2  + 
         2^{\frac{1}{2}} 3^{\frac{3}{4}})) $$
In a similar way we can get a polynomial $\bar{B}_{9}^{2}(z)$ with three critical points with critical value -1, corresponding to the black vertices in Fig.1b, $q=3$. The roots of $\frac{d\bar{B}_{9}^{2}(z)}{dz}$ with  $\zeta=-1$ are denoted by $a,b,0$, while $u$ has $\zeta=1$. We obtain in this case $u^{9}+18=0$, and  $a,b$ are the roots of $ z^2 - 3 u z +3 u^2=0$. We notice that $\bar{B}_{9}^{2}(z)$ has complex coefficients even if we take $u \in \bf{R}$. It has the same planar tree as $G_{3,3,3}(z)$ in Prop.2, but the polynomial derived from the classical Jacobi polynomial is simpler.
\par
The surfaces $J_{9}^{C}(x,y)+G^{2}_{9}(z)=0$, with $G^{2}_{9}(z)=(B_{9}^{2}(z)+1)/2$ and $G^{2}_{9}(z)=(\bar{B}_{9}^{2}(z)+1)/2$, have 110 and 127 cusps respectively.  The surfaces $J_{3q}^{C}(x,y)+( \bar{B}_{3q}^{2}(z)+1)/2=0$ have $\lfloor \frac{3q-1}{2}\rfloor-\lfloor \frac{3q}{3}\rfloor$ more cusps than those constructed in \cite{lab06} with $F_{d}(x,y)$. 
\par 
Other lower bounds can be improved by using $J_{m}^{C}(x,y)$. For $m=9,j=4$  we define $\bar{B}_{9}^{4}(z)$ in such a way that $\frac{d\bar{B}_{9}^{4}(z)}{dz}=(z-a)^{4}(z-b)^{4}$, with $\bar{B}_{9}^{4}(a)=1,
\bar{B}_{9}^{4}(b)=-1$. The Groebner basis has 2 elements and if we take the real root of $315+128 b^9=0$, then $a=-b$. In contrast to the cases studied above we find a solution for $\bar{B}_{9}^{4}(z)$ with real coefficients. The surface $J_{9}^{C}(x,y)+(\bar{B}_{9}^{4}(z)+1)/2=0$, where the normalized Belyi polynomial is
$$\frac{1}{80640}(40320 +  6^{7/9} 35^{8/9} 945 z - 6^{1/3} 35^{2/3}5040 
 z^3 + 6^{8/9} 35^{4/9} 3024 z^5 -  6^{4/9} 35^{2/9}5760
z^7 + 4480 z^9),$$
\par\noindent 
has 55 real singularities of type $A_{4}$ (Fig.2c).

\bigskip\par 
 \section{References}
 \bigskip\par

\end{document}